\documentclass{daj}
\usepackage{amsmath,amsfonts,amsbsy,amsgen,amscd,mathrsfs,amssymb,amsthm}
\usepackage{enumerate,mathtools}

\newtheorem{thm}{Theorem}[section]
\newtheorem{theorem}{Theorem}
\newtheorem{corollary}{Corollary}
\newtheorem{lem}[thm]{Lemma}

\newtheorem{prop}[thm]{Proposition}

\newtheorem{defn}[thm]{Definition}

\newtheorem{rem}{Remark}
\dajAUTHORdetails{%
  title = {The KKL inequality  and  Rademacher type 2}, 
  author = {Paata Ivanisvili and Yonathan Stone},
  plaintextauthor = {Paata Ivanisvili, Yonathan Stone},
  plaintexttitle = {The KKL inequality  and  Rademacher type 2}, 
  runningtitle = {}, 
  runningauthor = {},
  copyrightauthor = {},
  keywords = {Rademacher type, KKL inequality, influences; Boolean functions},
}   

\dajEDITORdetails{%
   year={2024},
   number={2},
   received={25 October 2022},   
   published={3 May 2024},  
   doi={10.19086/da.117131},       
}   

\begin{document}

\begin{frontmatter}[classification=text]

\title{The KKL inequality  and  Rademacher type 2} 

\author[pi]{Paata Ivanisvili\thanks{Supported in part by NSF grants DMS-2152346 and 
CAREER-DMS-2152401}}
\author[ybest]{Yonathan Stone}

\begin{abstract}
We show that a vector-valued Kahn--Kalai--Linial inequality holds in every Banach space of Rademacher type 2. We also show that  for any nondecreasing function $h\geq 0$ with $0<\int_{1}^{\infty}\frac{h(t)}{t^{2}}\mathrm{dt}<\infty$ we have the inequality
\begin{align*}
 \|f - \mathbb{E}f\|_2 \leq 12 \, T_{2}(X) \left(\int_{1}^{\infty}\frac{h(t)}{t^{2}} \mathrm{dt} \right)^{1/2} 
   \, \left(\sum_{j=1}^n \frac{\|D_j f\|^{2}_2}{h\left( \log \frac{\|D_j f\|_2}{\|D_j f\|_1} \right)}\right)^{1/2}
\end{align*}
for all $f :\{-1,1\}^{n} \to X$ and all $n\geq 1$, where $X$ is a normed space and $T_{2}(X)$ is the associated type 2 constant. 
\end{abstract}
\end{frontmatter}

\section{Introduction}
One of the central results in the analysis of Boolean functions is an estimate due to Kahn, Kalai, and Linial from 1988 \cite{KKLOG}. For $n\geq 1$, consider the discrete hypercube $\{-1,1\}^n$ equipped with the uniform probability measure. Set $\varepsilon=(\varepsilon_{1}, \ldots, \varepsilon_{n}) \in \{-1,1\}^{n}$.   For any \textit{Boolean} function $f :\{-1,1\}^{n} \to \{-1,1\}$ with  variance $\mathrm{Var}(f) = \mathbb{E} |f-\mathbb{E} f|^{2} =\delta > 0$ one has 
\begin{align}\label{kkl01}
\max_{1\leq j \leq n} \mathrm{Inf}_{j}(f) \geq C(\delta) \frac{\log(n)}{n},
\end{align}
where $\mathrm{Inf}_{j}(f) =  \mathbb{P}(f(\varepsilon) \neq  f(\varepsilon^{\oplus j}))$, $\varepsilon^{\oplus j} = (\varepsilon_{1}, \ldots, -\varepsilon_{j} , \ldots, \varepsilon_{n})$, and the positive constant $C(\delta)$ depends only on $\delta$. In fact one can take $C(\delta) = \delta/5$, see Corollary~\ref{kk09} in the Appendix. The quantity $\mathrm{Inf}_{j}(f)$, also called the influence of the Boolean function $f$ in the $j$'th coordinate,  measures how {\em influential} the $j$'th variable (or $j$'th voter) is towards the outcome of $f$. Clearly, if $\mathrm{Inf}_{j}(f)=0$, then $f$ does not depend on the
$j$'th variable at all.

Perhaps the most striking application of the inequality (\ref{kkl01}) is that for any monotone Boolean function $f$, that is, one for which $f(\varepsilon) \leq f(\varepsilon')$ whenever $\varepsilon_{j} \leq \varepsilon'_{j}$ for all $j=1,\ldots, n$, and for which $\mathbb{E} f \geq -0.99$, one can take $O\left(n/\log(n+1)\right)$ variables to be $1$, resulting in the average of $f$ with respect to the remaining variables being larger than $0.99$ (see Proposition 9.27 in \cite{Ryan}).  In social choice, this means that in a monotone election between two candidates, if the first candidate (corresponding to ``$-1$'') wins the election on average, then the second candidate can bribe $O(n/\log(n+1))$ voters to win the election on average.

Working with Boolean functions one is likely to fall into the trap of assuming that many of the phenomena observed in (\ref{kkl01}) are combinatorial in nature. Note that the classical Poincar\'e inequality 
\begin{align}\label{poincare}
\mathbb{E} |f - \mathbb{E} f|^{2} \leq \sum_{j=1}^{n} \mathbb{E} | D_{j} f|^{2},
\end{align}
for all $f :\{-1,1\}^{n} \to \mathbb{R}$, where $D_{j}f(\varepsilon) = \frac{f(\varepsilon) - f(\varepsilon^{\oplus j})}{2}$, 
implies a weak version of (\ref{kkl01}), i.e., for any Boolean $f :\{-1,1\}^{n} \to \{-1,1\}$ there exists $j \in \{1, \ldots, n\}$ such that $\mathrm{Inf}_{j}(f) = \mathbb{E} |D_{j} f|^{2} \geq \mathbb{E} |f-\mathbb{E}f|^{2}/n$.

By all accounts, Talagrand ~\cite{Tal1} was the first to discover a more extreme version of this inequality, which itself implies the existence of a universal constant $K \in (0, \infty)$ such that


\begin{align}\label{ienflo}
\mathbb{E} |f-\mathbb{E} f|^{2} \leq \frac{K}{\log \left(e /  \max_{k} \frac{\mathbb{E} |D_{k} f|}{\sqrt{\mathbb{E} |D_{k} f|^{2}}} \right)} \sum_{j=1}^{n} \mathbb{E} |D_{j} f|^{2}
\end{align}
 holds for all real-valued functions $f :\{-1,1\}^{n} \to \mathbb{R}$. We should point out that  (\ref{ienflo}) is the version of this estimate that is used in all applications we are aware of. It also improves the Poincar\'e inequality (\ref{poincare}) by a logarithmic factor, and can also be used to recover (\ref{kkl01}) when applied to Boolean functions. The inequality (\ref{ienflo}) for Boolean functions, although never stated in this form, follows from the arguments of the original paper of Kahn--Kalai--Linial~\cite{KKLOG} (see the Appendix for the full derivation), and for this reason we shall refer to (\ref{ienflo}) as the KKL inequality. 

The main advantage of the estimate (\ref{ienflo}) over the original Boolean corollary (\ref{kkl01}) is that it can be restated when the target space of the function is taken to be an arbitrary metric space.

\begin{defn}
We say a metric space $(X,d)$ is of  \textit{KKL type} if there exists a universal constant $T\in (0,\infty)$ such that the inequality 
\begin{equation}  
\mathbb{E}d(f(\varepsilon),f(\varepsilon'))^2 \leq \frac{T^{2}}{\log\left(e/\max\limits_k\frac{\mathbb{E}d(f(\varepsilon),f(\varepsilon^{\oplus k}))}{\sqrt{\mathbb{E}d(f(\varepsilon),f(\varepsilon^{\oplus k}))^2}}\right)}\, \sum_{j=1}^{n} \mathbb{E}d(f(\varepsilon),f(\varepsilon^{\oplus j}))^2 \label{MetricKKL}
\end{equation}
 holds for every $n\geq 1$, and for every function $f:\{-1,1\}^n \to X$, where $\varepsilon'$ is an independent copy of $\varepsilon$. The best such $T$  in (\ref{MetricKKL}) is denoted by $T_{\mathrm{KKL}}(X)$. 
\end{defn}

Thus an analogous statement to (\ref{kkl01})  can be formulated without relying on the linear structure of the target space $X$ one has when working with normed spaces.  The translation of phenomena involving a linear structure to fit spaces endowed only with metrics lies at the foundation of the Ribe program \cite{Naor1} which was inspired by Ribe's rigidity theorem~\cite{Rib1}, and was initiated by Bourgain \cite{Bo1}.

This paper seeks to provide a complete description of normed spaces which are of KKL type.  It is straightforward to see that a normed space  $(X, \|\cdot \|)$ that is of KKL type must also be of \textit{Rademacher type 2} (or just type 2), meaning that there exists a constant $T\in (0,\infty)$ such that for all $x_1,...,x_n \in X$ we have that
\begin{align}\label{type2}
\mathbb{E}\Big\| \sum_{j=1}^{n} \varepsilon_{j} x_{j} \Big\|^{2} \leq T^{2}  \sum_{j=1}^{n} \|x_j\|^2.
\end{align}
The best constant $T$ in (\ref{type2}) is denoted by $T_{2}(X)$.

A recent  paper \cite{IVHV} showed that  normed spaces of type 2  must satisfy the Poincar{\'e} inequality (\ref{poincare}), resolving a long-standing conjecture in Banach space theory due to Enflo. Soon after, Eskenazis and Cordero-Erausquin showed \cite{Esk1} that for type 2 spaces there was a variant of the KKL inequality (\ref{ienflo}) including an additional doubly logarithmic factor, which meant that (\ref{kkl01}) could not be recovered from this. 
An instructive example, which serves as a ``self-checker''  for whether any given new functional inequality bounding the variance of a function $f$ by its discrete derivatives can ever recover (\ref{kkl01}) is the {\em tribes} function $f_{\mathrm{tribe}}: \{-1,1\}^{n} \mapsto \{-1,1\}$, which has the properties: $\mathbb{E} ( \sum_{j=1}^{n} \mathrm{Inf}_{j}(f))^{\alpha} \asymp \log^{\alpha}(n)$ for any $\alpha>0$; $\max_{j=1,...,n} \mathrm{Inf}_{j}(f_{\mathrm{tribe}}) = \mathrm{Inf}_{1}(f_{\mathrm{tribe}}) \asymp \frac{\log(n)}{n}$; $\mathrm{Var}(f_{\mathrm{tribe}}) \asymp 1$ (see ~\cite{BEN}). Notice that for the tribes function the right hand side of (\ref{ienflo}) is of constant order as $n \to \infty$.

In this paper we show that KKL type and Rademacher type 2 coincide. 

\begin{theorem}\label{mth1} For any normed space $(X, \|\cdot \|)$, we have
\begin{align*}
 T_{2}(X)/2 \leq T_{\mathrm{KKL}}(X) \leq  2e \sqrt{2\pi}\, T_{2}(X).
\end{align*}

\end{theorem}

It will follow from the proof of Theorem~\ref{mth1} that if one manages to obtain bounds 
$\| D_{j} f\|_{2}\leq b_{j}, \|D_{j} f\|_{1} \leq a_{j}$ for some nonnegative numbers $a_{j}, b_{j}$ with $a_{j}\leq b_{j}$, $j=1, \ldots, n$, 
then the following inequality holds 
\begin{align}\label{gen1}
 \| f - \mathbb{E} f\|_{2} \leq \frac{2e\sqrt{2\pi} \, T_{2}(X)}{\log(e/\max_{j} (a_{j}/b_{j}))} \Big(\sum_{k=1}^{n} b_{k}^{2}\Big)^{1/2},
\end{align}
where $T_{2}(X)$ is the type 2 constant of the normed space $X$. 

The inequality, as written in (\ref{gen1}), applied to real-valued functions was critical (see Chapter 5 in \cite{SC}) to obtaining
 sublinear bounds on the variance in the first passage percolation model \cite{BKS}.

\vskip0.5cm 

One may wonder whether there is an analog of Theorem~\ref{mth1} for normed spaces of type $p$, for $p \in [1,2]$. Recall that the normed space $(X, \| \cdot \|)$ is of type $p$,  $p \in [1,2]$, if there exists a positive constant $T$ such that 
\begin{align}\label{radp}
\mathbb{E} \Big\| \sum_{j=1}^{n} \varepsilon_{j} x_{j} \Big\|^{p} \leq T_{p}^{p} \sum_{j=1}^{n} \|x_{j}\|^{p}
\end{align} 
for all $n\geq 1$ and all $x_{1}, \ldots, x_{n} \in X$. The best constant $T$ in (\ref{radp}) is denoted by $T_{p}(X)$. 

\begin{theorem}\label{mth2}
If a normed space $(X, \|\cdot \|)$ is of type $p$, then we have 
\begin{align} \label{obsh}
\|f-\mathbb{E} f\|_{p} \leq \frac{2e\sqrt{2\pi}\,  T_{p}(X)}{\sqrt{\log(e/\max (\frac{a_{j}}{b_{j}}))}} \left(\sum_{k=1}^{n}b_{j}^{p} \right)^{1/p}
\end{align}
 for any $n\geq 1$, any $f :\{-1,1\}^{n} \to X$, and any nonnegative numbers $\{a_{j}, b_{j}\}_{j=1}^{n}$, $b_{j}\geq a_{j}$ with $\| D_{j} f\|_{1}\leq a_{j}$, $\|D_{j} f\|_{p} \leq b_{j}$, $j=1, \ldots, n$. 
\end{theorem}

Our third theorem investigates Talagrand's inequality for normed spaces of type 2. Theorem 1 in \cite{Esk1} says that if $X$ is of type 2 then 
\begin{align}\label{dar1}
    \|f - \mathbb{E}f\|_2 \leq \frac{C}{\sqrt{\varepsilon}}\left(\sum_{j=1}^n \frac{\|D_j f\|^{2}_2}{ \log^{1-\varepsilon}\frac{\|D_j f\|_2}{\|D_j f\|_1}}\right)^{1/2}
\end{align}
holds for all $\varepsilon \in (0,1)$ and all $f : \{-1,1\}^{n} \mapsto X$, where $C$ depends only on the type 2 constant of $X$.  It is an open problem whether one can remove $\varepsilon$ in the statement of the inequality (\ref{dar1}), i.e., remove the constant $\frac{1}{\sqrt{\varepsilon}}$ and take $\varepsilon =0$ in the power of the logarithm. 

In this direction we make the following improvement.

\begin{theorem}\label{mth3}
Let $h \geq 0$ be nondecreasing with $0<\int_1^\infty h(t) \tfrac{dt}{t^2} < \infty$ and $X$ a normed space of type 2.  Then 
\begin{align}\label{g-Tal}
    \|f - \mathbb{E}f\|_2 \leq 12 \, T_{2}(X) \left(\int_{1}^{\infty}\frac{h(t)}{t^{2}} \mathrm{dt} \right)^{1/2} 
   \, \left(\sum_{j=1}^n \frac{\|D_j f\|^{2}_2}{h\left( \log \frac{\|D_j f\|_2}{\|D_j f\|_1} \right)}\right)^{1/2}
\end{align}
for all $f :\{-1,1\}^{n} \to X$ and all $n\geq 1$. 
\end{theorem}

The choice $h(t) = t^{1-\varepsilon}$ recovers (\ref{dar1}). One can consider  more sophisticated examples such as $h(t) = \frac{t}{\log^{1+\varepsilon}(2+t)}$ or $h(t) = \frac{t}{\log(2+t) (\log \log (10+t))^{1+\varepsilon}}$ for any $\varepsilon>0$. 

  It is worth noting that if a normed space $X$ satisfies (\ref{g-Tal}) for some $h$ as above, then by considering linear functions $f(\varepsilon) = \sum_i\varepsilon_ix_i$, we see that $X$ must be of type 2. 
  
  In Sections \ref{ss45}, \ref{ss21}, \ref{ss23} we present proofs of Theorems~\ref{mth1}, \ref{mth2}, \ref{mth3}, respectively. In Section~\ref{ss25} we present Propositions \ref{utv1} and \ref{utv2} showing that the main steps in the proofs of our theorems give sharp bounds.


\section{Background}

Consider the space of vector-valued functions $f:\{-1,1\}^n \rightarrow X$, defined on the discrete hypercube, for which $(X,\|\cdot\|)$ is a normed space.  Define the $L^p(X)$ norm on such functions as 
\[\|f\|_{p}  = \left(\mathbb{E}\|f\|^p \right)^{1/p} = \Big(\frac{1}{2^n}\sum\limits_{\varepsilon \in \{-1,1\}^n} \|f(\varepsilon)\|^p \Big)^{1/p}.\]
It is often helpful to think about $f$ in terms of its  \textit{Fourier-Walsh} expansion
\[f(\varepsilon) = \sum_{S \subset \{1,...,n\}} a_S\varepsilon^S,\]
where 
\[\varepsilon^S: (\varepsilon_1,...,\varepsilon_n) \longmapsto \prod\limits_{i \in S} \varepsilon_{i}, \quad \varepsilon^{\emptyset} \equiv 1,\]
are called \textit{Walsh functions}.  A  straightforward computation yields 
\[a_S = \mathbb{E}\varepsilon^{S}f(\varepsilon) \quad \text{ for all } \quad S\subseteq \{1,...,n\}.\]
Much of what we will prove in this paper relies on the rich theory of the \textit{heat semigroup} on the discrete hypercube.  With the discrete derivative operators $D_j$ defined as in the introduction, define the \textit{discrete Laplacian} $\Delta$ on functions $f:\{-1,1\}^n \to X$ as follows
\[\Delta f(\varepsilon)= - \sum_{j=1}^{n} D_{j} f(\varepsilon).\]
We then define the \textit{heat semigroup} $P_t = e^{\Delta t}$ as
\[P_{t}(f)(\varepsilon) = \sum_{S \subset \{1,...,n\}} a_Se^{-|S|t}\varepsilon^S \quad \text{for all}\quad t \geq 0.\]
An important property satisfied by the heat semigroup is \textit{hypercontractivity}: for all $p,q$ satisfying $1 < p \leq q < \infty$ and $e^{-2t} \leq \frac{p-1}{q-1}$, we have   
\begin{align}\label{hyp1}
\|P_t f\|_{q} \leq \|f\|_{p}.
\end{align}

One of the key ingredients in obtaining Theorem~\ref{mth1}  is going to be the following pointwise identity obtained in \cite{IVHV}.  For a normed space $(X,||\cdot||)$ and a function $f:\{-1,1\}^n \to X$, we have 
\begin{equation}
    -\frac{d}{dt}P_t f(\varepsilon) = \frac{1}{\sqrt{e^{2t} - 1}}\mathbb{E}_\xi \sum_{j=1}^n \delta_j(t) D_j f(\varepsilon\xi(t)), \label{heat}
\end{equation}
where $\varepsilon\xi(t) = (\varepsilon_1\xi_1(t),...,\varepsilon_n\xi_n(t))$, and the $\xi_i(t)$ are i.i.d. random variables with
\[\mathbb{P}\{\xi_i(t) = \pm 1\} = \frac{1\pm e^{-t}}{2},\]
and   $\delta_i = \frac{\xi_i(t) - \mathbb{E} \xi_{i}(t)}{\sqrt{\mathrm{Var}(\xi_{i}(t))}}$.  

\section{Proof of Theorem~\ref{mth1}}\label{ss45}
The lower bound $T_{2}(X)/2 \leq T_{\mathrm{KKL}}(X)$ follows by applying (\ref{MetricKKL}) to linear functions $f(\varepsilon)=\varepsilon_{1} x_{1}+\ldots+\varepsilon_{n} x_{n}$ and using the simple inequality $\mathbb{E}\| f-\mathbb{E}f \|^{2} \leq \mathbb{E} \| f(\varepsilon) - f(\varepsilon')\|^{2}$.

  To get the upper bound in Theorem~\ref{mth1} we proceed as follows. Without loss of generality we assume $\max_{1\leq k \leq n} \mathbb{E} \| D_{k} f\| >0$ otherwise $f$ is constant and hence there is nothing to verify.   By applying (\ref{heat}) to the function $P_t(f)$, and using the chain rule, we obtain 
   \[-\frac{d}{dt}P_{2t} f(\varepsilon) = \frac{2}{\sqrt{e^{2t} - 1}}\mathbb{E}_\xi \sum_{j=1}^n \delta_j(t) D_j P_tf(\varepsilon\xi(t)).\]
   We then proceed by integrating both sides with respect to $t$ and note that $P_{0}f = f$,  and $\lim_{t \to \infty}P_{t}f = \mathbb{E}f$.
   It follows that
   \begin{align}\label{chainr}
   f(\varepsilon) - \mathbb{E}f = 2\cdot\int_0^\infty \mathbb{E}_\xi \sum_{j=1}^n \delta_j(t) D_j P_tf(\varepsilon\xi(t))\frac{\mathrm{d}t}{\sqrt{e^{2t}-1}}.
   \end{align}
   Next, we take $L^2$ norms and apply both the triangle and Jensen's inequality, which yields
   \[\|f-\mathbb{E}f\|_2\leq 2\cdot\int_0^\infty \Big(\mathbb{E}_{\xi,\varepsilon}\Big\| \sum_{j=1}^n \delta_j(t) D_j P_tf(\varepsilon\cdot\xi(t))\Big\|^2\Big)^{1/2}\frac{\mathrm{d}t}{\sqrt{e^{2t}-1}}.\]
   Since $(\mathbb{E} \| f(\varepsilon) - f(\varepsilon')\|^{2})^{1/2} \leq 2  \| f-\mathbb{E} f\|_{2}$,  we will shift to exclusively obtaining upper bounds on the right hand side (RHS) of the inequality above. Notice that  $(\delta,\xi\varepsilon)$ has the same  distribution as $(\delta,\varepsilon)$.  We thus replace $\xi\varepsilon$ with $\varepsilon$ via a change of variables, i.e.
   \[\|f-\mathbb{E}f\|_2\leq \text{RHS} :=  2\cdot\int_0^\infty \Big(\mathbb{E}_{\xi,\varepsilon}\Big\| \sum_{j=1}^n \delta_j(t) D_j P_tf(\varepsilon)\Big\|^2\Big)^{1/2} \frac{\mathrm{d}t}{\sqrt{e^{2t}-1}}.\]
   In addition, a symmetrization argument using $\mathbb{E}\delta_i(t) = 0$ yields the following inequality
   \[\text{RHS} \leq  2\cdot\int_0^\infty \Big(\mathbb{E}_{\xi,\xi',\varepsilon,\varepsilon'}\Big\| \sum_{j=1}^n \varepsilon'_j(\delta_j(t)-\delta_j'(t)) D_j P_tf(\varepsilon)\Big\|^2\Big)^{1/2} \frac{\mathrm{d}t}{\sqrt{e^{2t}-1}},\]
   where $\delta_j'$ is a copy of $\delta_j$ depending on $\xi'$, which is an independent copy of $\xi$, and $\varepsilon' \in \{-1,1\}^n$ is an independent copy of $\varepsilon$.  Since our normed space $X$ is of type 2, applying the definition to the average over $\varepsilon'$ we obtain
   \begin{equation}
       \text{RHS} \leq 2T_{2}(X)\cdot\int_0^\infty \Big(\mathbb{E}_{\xi,\xi',\varepsilon} \sum_{j=1}^n (\delta_j(t)-\delta'_j(t))^2 \cdot \|D_j P_tf(\varepsilon)\|^2 \Big)^{1/2} \frac{\mathrm{d}t}{\sqrt{e^{2t}-1}}. \label{inter}
   \end{equation}
   Next, we have $\mathbb{E}_{\xi',\xi}(\delta_j(t) - \delta_j'(t))^2 = 2$,  which we can plug into (\ref{inter}) to obtain the following estimate
   \[\text{RHS} \leq 2^{3/2}\cdot T_{2}(X) \cdot\int_0^\infty \Big( \sum_{j=1}^n \|D_j P_tf\|^{2}_{2} \Big)^{1/2} \frac{\mathrm{d}t}{\sqrt{e^{2t}-1}}.\]
   Moreover, the commutativity $D_jP_t = P_tD_j$ along with hypercontractivity gives us
   \[\|D_j P_tf\|^{2}_{2} \leq \|D_j f\|_{1+e^{-2t}}^2.\]
   It thus follows that
   \[\text{RHS} \leq 2^{3/2}\cdot T_{2}(X) \cdot\int_0^\infty \Big( \sum_{j=1}^n \|D_j f\|_{1+e^{-2t}}^2 \Big)^{1/2} \frac{\mathrm{d}t}{\sqrt{e^{2t}-1}}.\]
   By H{\"o}lder's inequality, we have 
   \begin{align}\label{hohe}
    \|D_j f\|_{1+e^{-2t}}^2 \leq \|D_j f\|_{1}^{2\theta}\cdot \|D_j f\|_{2}^{2(1-\theta)} \leq \|D_j f\|_{2}^2 \cdot \left(\max_k\frac{\|D_k f\|_{1}}{\|D_k f\|_{2}}\right)^{2\theta},
    \end{align}
    where $\theta = \frac{1-e^{-2t}}{1+e^{-2t}}.$
    We have therefore obtained the following inequality
    \begin{equation}
        \|f-\mathbb{E}f\|_2\leq 2^{3/2}\cdot T_{2}(X) \cdot  \Big( \sum\limits_{j=1}^n\|D_jf\|_2^2\Big)^{1/2}\cdot \int_0^\infty \left(\max_k\frac{\|D_k f\|_1}{\|D_k f\|_2}\right)^{\frac{1-e^{-2t}}{1+e^{-2t}}}\frac{\mathrm{d}t}{\sqrt{e^{2t}-1}}. \label{pre-final}
    \end{equation}
    Let $a := \max_k\frac{\|D_k f\|_1}{\|D_k f\|_2} \in (0,1]$.  Since $\frac{1+e^{-t}}{1+e^{-2t}} \geq 1$ and $\frac{1}{\sqrt{1+e^{-t}}}\leq 1$, we have
    \[\int_0^\infty a^{\frac{1-e^{-2t}}{1+e^{-2t}}}\frac{\mathrm{d}t}{\sqrt{e^{2t}-1}} \leq \int_0^\infty a^{1-e^{-t}} \frac{e^{-t}\mathrm{d}t}{\sqrt{1-e^{-t}}} \leq e\cdot\int_0^\infty (a/e)^{1-e^{-t}} \frac{e^{-t}\mathrm{d}t}{\sqrt{1-e^{-t}}},\]
    Performing the substitution $s = \sqrt{(1-e^{-t})\log(e/a)}$ then gives us
    \[e\cdot\int_0^\infty (a/e)^{1-e^{-t}} \frac{e^{-t}\mathrm{d}t}{\sqrt{1-e^{-t}}} = \frac{2e}{\sqrt{\log(e/a)}}\int_0^{\sqrt{\log(e/a)}}e^{-s^2}\mathrm{d}s \leq \frac{e\sqrt{\pi}}{\sqrt{\log(e/a)}},\]
    whence the desired result follows immediately.
    \hfill \qedsymbol

\section{Proof of Theorem~\ref{mth2}}\label{ss21}

The proof of Theorem~\ref{mth2} proceeds in the same way as the proof of Theorem~\ref{mth1}, therefore,  to avoid repetitions we briefly sketch the argument. We take $L^{p}$ norms in both sides of (\ref{chainr}), and apply the triangle and Jensen's inequality, as well as the same symmetrization arguments.  By Cauchy--Schwarz, we have  $\mathbb{E} |\delta_{j}(t)-\delta'_{j}(t)|^{p} \leq 2^{p/2}$.  Then using hypercontractivity and  H\"older's inequality we obtain 
\begin{align*}
 \| D_{j} P_{t} f\|_{p}^{p} \leq \|D_{j} f\|_{1+(p-1)e^{-2t}}^{p} \leq b_{j}^{p} \left(\max_{k}\frac{a_{k}}{b_{k}}\right)^{p \theta}, \quad \text{where} \quad \theta = \frac{1-e^{-2t}}{1+(p-1)e^{-2t}}.
\end{align*}
 Finally note that for any $a \in (0,1]$  we have 
\begin{align*}
\int_{0}^{\infty} a^{\frac{1-e^{-2t}}{1+(p-1)e^{-2t}}} \frac{\mathrm{d}t}{\sqrt{e^{2t}-1}} \leq \int_{0}^{\infty} a^{1-e^{-t}} \frac{e^{-t}\mathrm{d}t}{\sqrt{1-e^{-t}}},
\end{align*}
after which the rest of the proof proceeds verbatim as in Section~\ref{ss45}.

Let us remark that the converse implication to Theorem~\ref{mth2} also holds true: if $(X, \| \cdot \|)$ is an arbitrary normed space, and  the inequality (\ref{obsh}) holds with $2e\sqrt{2\pi}\,  T_{p}(X)$ replaced by some universal constant, then the normed space $X$ must be of type  $p$. Indeed, the conclusion follows by choosing  $f = \sum_{j} \varepsilon_{j} x_{j}$ and  $a_{j}= b_{j} = \|x_j\|$.

\section{Proof of Theorem~\ref{mth3}} \label{ss23}
The first steps of the proof proceed along the same lines as in Theorem~\ref{mth1}, except that we do not take the maximum over $k$ in (\ref{hohe}).  We thus end up with the following inequality:
\begin{align}\label{TalStart}
    \|f-\mathbb{E}f\|_2\leq 2^{3/2}\cdot T_{2}(X) \cdot  \int_0^\infty \biggl(\sum\limits_{j=1}^n\|D_jf\|_2^2 \cdot  \left(\frac{\|D_j f\|_1^2}{\|D_j f\|_2^2}\right)^{\frac{1-e^{-2t}}{1+e^{-2t}}}\biggr)^{1/2}\frac{\mathrm{d}t}{\sqrt{e^{2t}-1}}. 
\end{align}
We then note that since $\frac{1-e^{-2t}}{1+e^{-2t}} \geq 1-e^{-t}$ and $\frac{\|D_j f\|_1}{\|D_j f\|_2} \leq 1$, as well as the fact that $\frac{1}{\sqrt{e^{2t} - 1}} \leq \frac{e^{-t}}{\sqrt{1-e^{-t}}}$, we can bound the right hand side in  (\ref{TalStart}) by
\[ \int_0^\infty \biggl(\sum\limits_{j=1}^n\|D_jf\|_2^2 \cdot  \left(\frac{\|D_j f\|_1}{\|D_j f\|_2}\right)^{1-e^{-t}}\biggr)^{1/2}\frac{e^{-t}\mathrm{d}t}{\sqrt{1-e^{-t}}}.\]
Using this estimate, we can perform the substitution $s = \sqrt{1-e^{-t}}$, which allows us to rewrite the integral above as
\[\frac{1}{2}\int_0^1 \biggl(\sum\limits_{j=1}^n\|D_jf\|_2^2 \cdot  \left(\frac{\|D_j f\|_1}{\|D_j f\|_2}\right)^{s^2}\biggr)^{1/2}\mathrm{d}s.\]
Let us denote $g(\sqrt{s})=h(s)$. Note that our theorem is satisfied if we can show 
\begin{align} \label{middle1}
    \int_0^1 \biggl(\sum\limits_j\|D_jf\|_2^2 \cdot  \left(\frac{\|D_j f\|_1}{\|D_j f\|_2}\right)^{s^2}\biggr)^{1/2}\mathrm{d}s \leq 12\sqrt{2} \left( \int_{1}^{\infty} \frac{g(s)}{s^{3}} \mathrm{ds} \right)^{1/2}\,   \biggl(\sum_{j=1}^n \frac{\|D_j f\|^2}{g\left(\log^{1/2} \tfrac{\|D_jf\|_2}{\|D_jf\|_1}.
    \right)}\biggr)^{1/2}.
\end{align}
 Next, denote   $c_j = \frac{\|D_jf\|_2^2}{g\left(\log^{1/2} \tfrac{\|D_jf\|_2}{\|D_jf\|_1}
    \right)},$ and $X_j = \tfrac{\|D_jf\|_1}{\|D_jf\|_2}$ for all $j=1,\ldots, n$. Then, dividing both sides of (\ref{middle1}) by $\Big(\sum\limits_{j=1}^n c_j\Big)^{1/2}$ gives us the following:
    \[\int_0^1 \Biggl(\sum\limits_{j=1}^n\Biggl(\frac{c_j}{\sum\limits_{k=1}^n c_k}\Biggr) \cdot X_j^{s^2} g(\log^{1/2}(1/X_j)) \Biggr)^{1/2}\mathrm{d}s \leq 12\sqrt{2} \, \left(\int_{1}^{\infty} \frac{g(s)}{s^{3}} \mathrm{ds}\right)^{1/2}.\]
    Viewing the $\frac{c_j}{\sum\limits_{k=1}^n c_k}$ as probabilities, we can obtain the result above, and thus prove Theorem~\ref{mth3} immediately by proving the following:
    \begin{lem}\label{texnika}
    For all $g \geq 0$ with $\int_{1}^{\infty} \frac{g(t)}{t^{3}}\mathrm{dt}<\infty$, we have that 
    \begin{align}\label{texn1}
        \int_0^1 \left(\mathbb{E} \: X^{s^2} g(\log^{1/2}(1/X))\right)^{1/2} \mathrm{d}s &\leq (\sqrt{2}+8\sqrt{\pi})\,  \left( \int_{1}^{\infty} \frac{g(s)}{s^{3}} \mathrm{ds}\right)^{1/2}
    \end{align}
    for all random variables $0 \leq X \leq 1$.
    \end{lem}
To prove this, rewrite $X = e^{-Y^2}$, where $Y$ is a non-negative random variable.  We thus have that
\begin{align*}
    \int_0^1 \left(\mathbb{E}\: X^{s^2} g(\log^{1/2}(1/X))\right)^{1/2} \mathrm{d}s = \int_0^1 \left(\mathbb{E} \: e^{-(sY)^2} g(Y)\right)^{1/2} \mathrm{d}s.
\end{align*}
Now write $p_k = \mathbb{P}( 2^{k}\leq Y< 2^{k+1})$ and use the fact that $t \mapsto e^{-t^2}$ is decreasing and $t \mapsto g(t)$ is nondecreasing for $t\geq 0$, we get
\[\int_0^1 \left(\mathbb{E} \: e^{-(sY)^2} g(Y)\right)^{1/2} \mathrm{d}s \leq \int_0^1 \Big(g(1)+\sum\limits_{k=0}^\infty e^{-(s2^k)^2} g(2^{k+1})p_k\Big)^{1/2} \mathrm{d}s.\]
Next, we partition the interval we are integrating over into intervals of the form $(2^{-j}, 2^{-j+1}]$ to get the following estimate
\begin{align*} 
\int_0^1 \Big(g(1)+\sum_{k=0}^\infty e^{-(s2^k)^2} g(2^{k+1})p_k\Big)^{1/2} \mathrm{d}s \leq \sum\limits_{j=1}^\infty \Big(g(1)+\sum\limits_{k=0}^\infty e^{-2^{2(k-j)}} g(2^{k+1})p_k\Big)^{1/2} 2^{-j}.
\end{align*}
Next, using subadditivity of the map $t \mapsto \sqrt{t}$ and interchanging the sums we obtain
\begin{align}\label{middle2}
    \sum\limits_{j=1}^\infty \Big(g(1)+\sum\limits_{k=0}^\infty e^{-2^{2(k-j)}} g(2^{k+1})p_k\Big)^{1/2} 2^{-j} \leq     \sqrt{g(1)}+\sum\limits_{k=0}^\infty \sqrt{g(2^{k+1})}\sqrt{p_k} \sum\limits_{j=1}^\infty e^{-\frac{2^{2(k-j)}}{2}} 2^{-j}.
\end{align}
Next, notice that
\[\sum\limits_{j=1}^\infty e^{-\frac{2^{2(k-j)}}{2}} 2^{-j} \leq \sum\limits_{j=1}^\infty \int_{2^{-j}}^{2^{-j+1}} e^{-\frac{(2^k x)^2}{8}} \mathrm{d}x = \int_0^1 e^{\frac{-(2^k x)^2}{8}}\mathrm{d}x = 2^{-k}\int_0^{2^k} e^{-y^2/8} \mathrm{dy}, \]
where the estimate $\int_0^{2^k} e^{-y^2/8} \mathrm{dy} \leq \sqrt{2\pi}$ gives us that 
\[\sum\limits_{j=1}^\infty e^{-\frac{2^{2(k-j)}}{2}} 2^{-j} \leq \sqrt{2\pi}\cdot 2^{-k},\]
for all $k$.  Using this estimate on (\ref{middle2}) gives us
\begin{align*}
&\sqrt{g(1)}+\sum\limits_{k=0}^\infty \sqrt{g(2^{k+1})}\sqrt{p_k} \sum\limits_{j=1}^\infty e^{-\frac{2^{2(k-j)}}{2}} 2^{-j} \leq \\
&\sqrt{g(1)}+\sqrt{2\pi}\cdot \sum\limits_{k=0}^\infty \sqrt{g(2^{k+1})}2^{-k}\sqrt{p_k}.
\end{align*}
However, notice that 
\[\sum_{k=0}^\infty g(2^{k+1})2^{-2k} \leq 8\sum_{k=0}^{\infty} \int_{2^{k}}^{2^{k+1}} \frac{g(2t)}{t^3}  \mathrm{dt} =32 \int_{2}^{\infty} \frac{g(s)}{s^{3}} \mathrm{ds}, \]
and $g(1) \leq 2 \int_{1}^{\infty}\frac{g(s)}{s^{3}}ds$. We thus have by Cauchy-Schwarz that
\[\sqrt{g(1)}+\sqrt{2\pi}\cdot \sum\limits_{k=0}^\infty \sqrt{g(2^{k+1})}2^{-k}\sqrt{p_k} \leq (\sqrt{2}+8\sqrt{\pi}) \left( \int_{1}^{\infty} \frac{g(s)}{s^{3}} \mathrm{ds}\right)^{1/2}.\]
which completes the proof of the lemma, and thus also the theorem. \hfill\qedsymbol

\section{Concluding remarks}\label{ss25}
One may wonder how sharp the bound obtained in Lemma~\ref{texnika} is. We can show that the inequality (\ref{texn1}) is sharp up to a multiplicative constant. 
\begin{prop}\label{utv1}
For any $g\geq 0$ with $\int_{1}^{\infty} \frac{g(t)}{t^{3}}\mathrm{dt}<\infty$ there exists a random variable $X$, $0\leq X\leq 1$, such that 
\begin{align*}
\int_{0}^{1} \left( \mathbb{E} X^{s^{2}} g(\log^{1/2}(1/X))\right)^{1/2}\mathrm{ds} \geq \frac{e^{-8}}{2}\left(\int_{1}^{\infty} \frac{g(s)}{s^{3}} \mathrm{ds}\right).
\end{align*}
\end{prop}
\begin{proof}
Let $X=e^{-Y^{2}}$ and choose $Y\geq 0$ so that   $p_{k} = \mathbb{P}(2^{k} \leq Y < 2^{k+1}) = \frac{g(2^{k}) 2^{-2k}}{\sum_{j=1}^{\infty} g(2^{j}) 2^{-2j}}$ for all $k\geq 1$. We have 
\begin{align*}
\int_{0}^{1} \left( \mathbb{E} X^{s^{2}} g(\log^{1/2}(1/X))\right)^{1/2}\mathrm{ds}&=\int_{0}^{1}\left(\mathbb{E} e^{-(sY)^{2}} g(Y)\right)^{1/2}\mathrm{ds} \geq \\
\sum_{j=1}^{\infty} \int_{2^{-j}}^{2^{-j+1}} \Big(\sum_{k=0}^{\infty}e^{-(s2^{k+1})^{2}}g(2^{k})p_{k} \Big)^{1/2}ds &\geq \sum_{j=1}^{\infty} \Big(\sum_{k=0}^{\infty} e^{-(2^{-j+k+2})^{2}} g(2^{k})p_{k} \Big)^{1/2}2^{-j} \geq \\
\sum_{j=1}^{\infty}e^{-8} \sqrt{g(2^{j}) p_{j}} \, 2^{-j} &= e^{-8} \Big(\sum_{j=1}^{\infty} g(2^{j}) 2^{-2j} \Big)^{1/2},
\end{align*}
where we established the final inequality by considering only the $k=j$ terms.  On the other hand 
\begin{align*}
\sum_{j=1}^{\infty} g(2^{j}) 2^{-2j} \geq  \sum_{j=1}^{\infty} \int_{2^{j}}^{2^{j+1}} \frac{g(t/2)}{t^{3}} \mathrm{dt} = \frac{1}{4}\int_{1}^{\infty} \frac{g(s)}{s^{3}}\mathrm{ds}.
\end{align*}
This finishes the proof of the proposition. 
\end{proof}

\vskip1cm 

One may also wonder whether we can prove the following inequality
\begin{align}\label{ale1}
\int_0^\infty \Biggl(\sum\limits_{j=1}^{n} \|D_jf\|_2^2 \cdot  \left(\frac{\|D_j f\|_1^2}{\|D_j f\|_2^2}\right)^{\frac{1-e^{-2t}}{1+e^{-2t}}}\Biggr)^{1/2}\frac{\mathrm{d}t}{\sqrt{e^{2t}-1}} \leq C \Biggl(\sum_{j=1}^n \frac{\|D_j f\|^{2}_2}{ \log \frac{\|D_j f\|_2}{\|D_j f\|_1}}\Biggr)^{1/2}
\end{align}
with some universal finite  constant $C>0$. Note that (\ref{ale1}) combined with (\ref{TalStart}) would solve the open problem regarding setting $\varepsilon=0$ in (\ref{dar1}).  If the inequality (\ref{ale1}) does not hold then this would mean that one needs to come up with a different approach to prove or disprove the inequality (\ref{dar1}) without $\varepsilon$. 

It was noted by the referee in \cite{Esk1} that if we treat $b_{j} = \| D_{j} f\|_{2}^{2}$, $ a_{j} = \|D_{j} f\|_{1}^{2}$ as arbitrary nonnegative numbers satisfying $b_{j}\geq a_{j} \geq 0$, then the inequality (\ref{ale1}) does not hold in such generality. However, in general the numbers $b_{j}$ and $a_{j}$ may not be arbitrary, i.e., there could be some relations between $\| D_{j} f\|_{1}$ and $\| D_{j} f\|_{2}$. For instance, for a Boolean function $f :\{-1,1\}^{n} \to \{-1,1\}$ we have $b_{j}=\|D_{j} f\|_{2}^{2} = \| D_{j} f\|_{1} = \sqrt{a_{j}}$ for all $j=1, \ldots, n$. Therefore, it is interesting to ask whether the inequality (\ref{ale1}) holds in the ``Boolean case'':

\begin{align}\label{bul1}
\int_0^\infty \biggl(\sum\limits_{j=1}^{n}  b_{j}^{1+\frac{1-e^{-2t}}{1+e^{-2t}}}\biggr)^{1/2}\frac{\mathrm{d}t}{\sqrt{e^{2t}-1}}
 \leq C' \biggl(\sum_{j=1}^n \frac{b_{j}}{ \log (1/b_{j})}\biggr)^{1/2}
\end{align} 
for all $n\geq 1$ and all $b_{j} \in [0,1]$.  It turns out that:

\begin{prop}\label{utv2}
There is no finite universal constant $C'>0$ for which the inequality (\ref{bul1}) holds for all $b_{1}, \ldots, b_{n}\in [0,1]$ and all $n\geq 1$.  
\end{prop}

\begin{proof}
 Towards a contradiction, assume that (\ref{bul1}) holds true. We have $b_{j}^{1+\frac{1-e^{-2t}}{1+e^{-2t}}} \geq b_{j}^{1+2(1-e^{-t})}$. Performing the substitution $s=\sqrt{1-e^{-t}}$ the inequality (\ref{bul1}) implies 
\begin{align}
\int_{0}^{1}\biggl( \sum_{j=1}^{n}b_{j}^{1+2s^{2}}\biggr)^{1/2} \mathrm{ds}  \leq C \biggl(\sum_{j=1}^n \frac{b_{j}}{ \log (1/b_{j})}\biggr)^{1/2}.
\end{align}
It follows from homogeneity that for all Borel measurable random variable $Y\geq 0$ we have 
\begin{align}\label{red01}
\int_{0}^{1}\left(\mathbb{E} e^{-Y^{2}-2(Ys)^{2}}\right)^{1/2}\mathrm{ds}\leq C \biggl( \mathbb{E}\frac{e^{-Y^{2}}}{Y^{2}}\biggr)^{1/2}.
\end{align}
Next, choose $Y$ so that it takes only the values $2^{k}$, $k \in \mathbb{Z}$ with probabilities $p_{k}$, i.e., $\sum_{k \in \mathbb{Z}}p_{k}=1$. Also set $p_{k}=0$ for $k\leq 0$. Then the right hand side in (\ref{red01}) takes the form
\begin{align*}
C\, \Big( \sum_{k \geq 1} 2^{-2k} e^{-2^{2k}}p_{k} \Big)^{1/2}.
\end{align*}
The left hand side in (\ref{red01}) we can lower bound as 
\begin{align*}
\sum_{j=1}^{\infty} \int_{2^{-j}}^{2^{-j+1}} \Big( \sum_{k \geq 1} e^{-2^{2k}-2(2^{k}s)^{2}}p_{k}\Big)^{1/2} \mathrm{ds} &\geq \\
\sum_{j=1}^{\infty} \Big( \sum_{k \geq 1} e^{-2^{2k}-2(2^{k-j+1})^{2}}p_{k}\Big)^{1/2} 2^{-j} &\geq e^{-4}\sum_{j\geq 1} \Big( e^{-2^{2j}} p_{j}\Big)^{1/2} \, 2^{-j}.
\end{align*}
Denoting $q_{j} =\left( e^{-2^{2j}} p_{j}\right)^{1/2} \, 2^{-j}$, $j\geq 1$, we see that (\ref{red01}) implies the inequality $\sum_{j\geq 1} q_{j}\leq e^{4}C (\sum_{j\geq 1} q_{j}^{2})^{1/2}$, which by homogeneity (and slightly abusing notation) must hold  for all $q_{j} \geq 0$, resulting in the desired contradiction.

\end{proof}

\appendix 
\section*{Appendix}

The improved Poincar\'e inequality  (\ref{ienflo}) for Boolean functions follows from the arguments of Kahn--Kalai--Linial~\cite{KKLOG}. 
\begin{theorem}[Kahn--Kalai--Linial]\label{vaxr}
We have
\begin{align}\label{kk03}
\mathbb{E} |f-\mathbb{E}f |^{2} \leq \frac{4}{\log \left( e/ \max_{k} \mathrm{Inf}_{k}(f)\right)} \sum_{j=1}^{n} \mathrm{Inf}_{j}(f)
\end{align}
for all $n\geq 1$ and any $f :\{-1,1\}^{n} \to \{-1,1\}$.
\end{theorem}

\begin{proof}
Assume $\mathbb{E} |f-\mathbb{E}f |^{2}>0$ (and hence $\max_{k} \mathrm{Inf}_{k}(f) > 0$) since otherwise there is nothing to prove. 
Let $f(\varepsilon) = \sum_{S \subset \{1,\ldots, n\}}\widehat{f}(S)  \varepsilon^{S}$. We start from the identity 
\begin{align}\label{kk1}
\sum_{S \subset \{1,\ldots, n\}} |S|\widehat{f}(S)^{2} = \sum_{j=1}^{n}\mathbb{E}|D_{j} f|^{2}.
\end{align}
The idea is to apply the identity (\ref{kk1}) to $P_{t}f$ instead of $f$ and integrate with respect to the measure $e^{-t}dt$ over the ray $[0, \infty)$. Indeed, note that $\widehat{P_{t}f}(S) = e^{-t|S|}\widehat{f}(S)$. Therefore (\ref{kk1}) gives
\begin{align}\label{kk02}
\sum_{S \subset \{1, \ldots, n\}} e^{-2t|S|}|S| a_{S}^{2} \leq \sum_{j=1}^{n} \mathbb{E} |P_{t} D_{j} f|^{2} \stackrel{(\ref{hyp1})}{\leq} \sum_{j=1}^{n} (\mathrm{Inf}_{j}(f))^{\frac{2}{1+e^{-2t}}},
\end{align}
where we also used the fact that $D_{j}f \in \{-1,0,1\}$. If we let $s :=e^{-2t}$ and integrate (\ref{kk02}) in $s$ over $(0,1)$ we obtain 
\begin{align*}
\sum_{S\neq \emptyset } \frac{|S|}{|S|+1}a_{S}^{2} \leq \int_{0}^{1} \sum_{j=1}^{n} \Big(\mathrm{Inf}_{j}(f)\Big)^{\frac{2}{1+s}}ds \leq  \Big( \sum_{j=1}^{n} \mathrm{Inf}_{j}(f)\Big) \, \int_{0}^{1} \Big(\max_{k} \mathrm{Inf}_{k}(f)\Big)^{\frac{1-s}{1+s}}ds. 
\end{align*}
On the one hand we have 
\begin{align*}
\frac{1}{2}\mathbb{E}|f-\mathbb{E}f|^{2} = \frac{1}{2}\sum_{S\neq \emptyset} a_{S}^{2} \leq \sum_{S\neq \emptyset } \frac{|S|}{|S|+1}a_{S}^{2}.
\end{align*}
Next, letting $a :=\max_{k} \mathrm{Inf}_{k}(f) \in (0,1]$, we obtain
\begin{align*}
\int_{0}^{1} a^{\frac{1-s}{1+s}}ds &\stackrel{s=1-x}{=} \int_{0}^{1}a^{\frac{x}{2-x}}dx \stackrel{a\leq 1}{\leq} \int_{0}^{1}a^{\frac{x}{2}}dx \leq 
\sqrt{e} \int_{0}^{1}(a/e)^{\frac{x}{2}}dx \stackrel{-x\log(a/e)=y}{=}\frac{1}{-\log(a/e)}\int_{0}^{-\log(a/e)}e^{-y/2}dy \\
&\leq \frac{2}{-\log(a/e)}.
\end{align*}
This completes the proof of (\ref{kk03}).
\end{proof}
\begin{rem}
It follows from the proof above that $\int_{0}^{1} \left(\mathrm{Inf}_{j}(f)\right)^{\frac{2}{1+s}}ds \leq  \frac{2 \mathrm{Inf}_{j}(f)}{-\log(\mathrm{Inf}_{j}(f)/e)}$ for all $j=1,\ldots, n$. Therefore, if we do not take the maximum in the proof of Theorem~\ref{vaxr} we would obtain 
$$
\mathbb{E} |f-\mathbb{E}f |^{2} \leq 4 \sum_{j=1}^{n}  \frac{\mathrm{Inf}_{j}(f)}{\log \left( e/ \mathrm{Inf}_{j}(f)\right)},
$$
which is Talagrand's inequality for Boolean functions. 
\end{rem}

\begin{corollary}\label{kk09}
For any $n\geq 1$ and any $f :\{-1,1\}^{n} \to \{-1,1\}$ we have 
\begin{align*}
\max_{k} \mathrm{Inf}_{1\leq k \leq n}(f) \geq \frac{1}{5}\mathrm{Var}(f) \frac{\log(n)}{n},
\end{align*}
where $\mathrm{Var}(f) = \mathbb{E} |f-\mathbb{E}f|^{2}$. 
\end{corollary}
\begin{proof}
If $n=1$ or $\mathrm{Var}(f)=0$ there is nothing to prove. Henceforth, we assume $n\geq 2$ and $\mathrm{Var}(f) = 1-(\mathbb{E}f)^{2} =:\delta \in (0,1]$. Let $t := \max_{k} \mathrm{Inf}_{1\leq k \leq n}(f)>0$. By (\ref{kk03}) we have $\delta \leq \frac{4 n t}{\log(e/t)}$. 
Since $\frac{d}{ds} \frac{s}{\log(e/s)}=\frac{2-\log(s)}{(1-\log(s))^{2}}>0$ the map $s \mapsto \frac{s}{\log(e/s)}$ is increasing on $(0,1]$. Assume the contrary, i.e., $t<\frac{\delta \log(n)}{5 n}$. To get a contradiction it suffices to show
\begin{align}\label{ggg}
\delta > \frac{4}{5} \frac{\delta \log(n)}{\ln(5en/(\delta\log(n)))}.
\end{align} 
The inequality (\ref{ggg}) is the same as  $5en^{1/5}> \delta \log(n)$ for all $n\geq 2$. On the other hand, for all $n\geq 1$ we have  $\frac{d}{dn}(5en^{1/5} - \log(n)) = \frac{en^{1/5}-1}{n}>0$ (and $5en^{1/5}> \delta \log(n)$ holds at $n=1$), hence the inequality (\ref{ggg}) is proved. 
\end{proof}

\section*{Acknowledgments} 
The authors are grateful to an anonymous referee for helpful suggestions and comments that improved the writing of the paper. The authors would like to thank Alexandros Eskenazis for valuable comments and remarks. 


\bibliographystyle{amsplain}


\begin{dajauthors}
\begin{authorinfo}[pi]
  Paata Ivanisvili\\
 Department of Mathematics, University of California, Irvine
\\
  Irvine, CA, USA\\
  pivanisv@uci.edu \\
\end{authorinfo}
\begin{authorinfo}[ybes]
  Yonathan Stone\\
  Department of Mathematics, University of California, Irvine\\
  Irvine, CA, USA\\
  ystone@uci.edu
\end{authorinfo}
\end{dajauthors}

\end{document}